\documentclass[12pt,draft,onecolumn]{IEEEtran}

\usepackage{amsmath}
\usepackage{amsthm}
\usepackage{amsfonts}
\usepackage{amssymb}
\usepackage{mathrsfs}
\usepackage{soul}
\usepackage{url}
\usepackage[usenames,dvipsnames]{xcolor}
\usepackage{enumerate}
\usepackage{tikz}
\usetikzlibrary{matrix}

\usepackage{datetime}

\newcommand{\R}{\mathbb{R}}
\newcommand{\C}{\mathbb{C}}

\newcommand{\rank}{\mathrm{rank}}
\newcommand{\rk}{\mathrm{rk}}
\newcommand{\diag}{\mathrm{diag}}

\newcommand{\im}{\mathrm{im}}

\renewcommand{\cos}{\mathrm{c}}
\renewcommand{\sin}{\mathrm{s}}

\theoremstyle{definition}

\newtheorem{Theorem}{Theorem}

\newtheorem{Definition}{Definition}

\definecolor{Royalblue}{cmyk}{1,0.30,0.2,0.2}

\newcommand{\tp}{^{\top}}
\newcommand{\mtp}{^{-\top}}
\newcommand{\inv}{^{-1}}

\newcommand{\beq}{\begin{equation}}
\newcommand{\eeq}{\end{equation}}
\newcommand{\bea}{\begin{eqnarray}}
\newcommand{\eea}{\end{eqnarray}}

\newcommand{\nn}{\nonumber}

\def\bmat{\left[ \begin{array}}
\def\emat{\end{array} \right]}

\author{Giacomo Baggio, Augusto Ferrante%
\thanks{Giacomo Baggio is with the Dipartimento di Ingegneria dell'€™Informazione, Universit\`a di Padova, via Gradenigo, 6/B€" I-35131 Padova, Italy. E-mail: giacomo.baggio@studenti.unipd.it.
Augusto Ferrante is with the Dipartimento di Ingegneria dell€'Informazione, Universit\`a di Padova, via Gradenigo, 6/B€" I-35131 Padova, Italy. E-mail: augusto@dei.unipd.it.}}

\begin{document}

\title{Parametrization of Minimal Spectral Factors of Discrete-Time Rational Spectral Densities}
\maketitle


\begin{abstract} 
In this paper, the problem of providing a complete parametrization of the minimal spectral factors of a discrete-time rational spectral density is considered. The desired parametrization, given in terms of the all-pass divisors of  a certain all-pass function, is established in the most general setting: after several partial results, mostly in the continuous-time case, this is indeed the first complete parametrization obtained without resorting to  any facilitating assumption.  This result provides a positive answer to a conjecture raised in \cite{BF-2016}.
\end{abstract}


\section{Introduction and problem definition} 

Spectral factorization problem is ubiquitous in systems and control theory. Some of its applications can be found in optimal estimation and filtering \cite{AM-79}, stochastic realization theory \cite{LP-15,Lindquist-P-85-siam, Lindquist-P-91-jmsec}, robust and optimal control theory \cite{ZDG-96,Willems-1971}, to cite just a few. Classical methodologies for the solution of the spectral factorization problem were introduced by Kolmogorov and Wiener in the early 40's for the case of scalar spectral densities and generalized by Youla \cite{Youla-1961} to the  multivariate case. Many works have been focused on developing numerically reliable algorithms for the computation of various important spectral factors, see e.g. \cite{C-93,OA-11}.

In \cite{BF-2016}, a general result on discrete-time spectral factorization was established and two conjectures were left to further investigation:
one was answered in the affirmative in \cite{BF-min-2016}. The other concerns the parametrization of the set of minimal spectral factors (i.e. with minimal McMillan degree) and it will be proven true in the present paper.
The problem of parametrizing the set of minimal spectral factors has a long history, see \cite{FP-82,Furhmann-95,FMP-93,P-93,F-94,Picci-P-94,F-97,FG-98,P-Ran-01,P-Ran-02,P-Ran-02b,F-05}, to cite just a few contributions. In fact, as pointed out in \cite{Lindquist-P-79-siam,Lindquist-P-91-jmsec}, this is a fundamental step in {\em stochastic realization theory}. A stochastic realization is
a representation of a second-order discrete-time purely nondeterministic  stationary process $\{y(t)\}$ as the output of a linear state-space model driven by white noise. Up to uninteresting changes of basis, stochastic realizations of minimal complexity are in one-to-one correspondence with minimal spectral factors of the spectral density $\Phi(z)$ of $\{y(t)\}$. For this reason, the problem of parametrizing the minimal spectral factors of $\Phi(z)$ is crucial for the analysis and synthesis of different models of a given stochastic process. 
However, to the best of our knowledge, all the available parametrizations rely on some restrictive assumptions on the spectral density and a general result is still missing.
This paper is an attempt to fill this gap: we provide a parametrization of the set of minimal spectral factors of a discrete-time rational spectral density in terms of the all-pass divisors of an all-pass function, which we name {\em conjugate phase function}. In doing so, we do not require any assumption on the considered spectral density. In particular, our result applies to spectral densities that are rank-deficient and/or possess zeros/poles on the unit circle and/or are improper.
The basis of our parametrization is the conjugate phase function, an all-pass function that can be explicitly computed from the minimum-phase spectral factor and the maximum phase unstable spectral factor: these two ``extremal spectral factors" can, in turn, be explicitly calculated as discussed in \cite{BF-2016} for the input-output 
representation and in the works by Oar\u{a} and co-workers \cite{O-05,OM-13} for the state-space representation.
Therefore, thanks to these contributions, our abstract theoretical parametrization  result may indeed be used to
explicitly provide all the minimal spectral factors of a given spectral density and hence all the minimal representations  of the corresponding process $\{y(t)\}$.

{\noindent \bfseries Paper structure.} The paper is organized as follows. In Section II, we collect some preliminary definitions and results on the parametrization of minimal spectral factors. In Section III, we state and prove the main result of the paper. In Section IV, we illustrate our main result by means of a numerical example. Finally, in Section V, we draw some conclusive remarks.

{\noindent \bfseries Notation.}  We denote by $\R$, $\C$, $\R^{m\times n}$, $\R[z]^{m\times n}$, and $\R(z)^{m\times n}$ the set of real numbers, complex numbers, $m\times n$ real matrices, real polynomial $m\times n$ matrices, and $m\times n$ real rational matrix-valued functions, respectively. Moreover, $\overline{\C}:=\C\cup\{\infty\}$. Given $G\in\R^{m\times n}$, we denote by $G\tp$ its transpose, by $G^+$ its Moore-Penrose pseudo-inverse,  by $G^{-L}$ its left-inverse, by $\ker(G)$ its kernel, and by $\im(G)$ its image. For $G\in\R^{m\times m}$,  $\sigma(G)$ denote the spectrum of $G$ and $G\mtp:=[G\tp]\inv$.
Given $G(z)\in\R(z)^{m\times n}$, we let $G^*(z):=G\tp(1/z)$. If $G(z)\in\R(z)^{m\times m}$, we let $G^{-*}(z):=[G^{*}(z)]\inv$. 
 Let $G(z)\in\R(z)^{n\times m}$ and let $\alpha\in\overline{\C}$, we denote by $\delta(G;\alpha)$ the degree of the pole of $G(z)$ at $\alpha$ , with the convention that $\delta(G;\alpha)=0$ if $\alpha\in\overline{\C}$ is not a pole of $G(z)$. We recall that the normal rank of $G(z)$, denoted by $\rk(G)$, is defined as the rank almost everywhere in $\C$. 
Finally, we denote by $\delta_{M}(G)$ the McMillan degree of $G(z)$ and we recall that the latter is equal to $\delta_{M}(G)=\sum_{i} \delta(G;\alpha_{i})$, where $\{\alpha_{i}\}\subset\overline{\C}$ is the set of poles of $G(z)$, see e.g. \cite[pag. 466]{K-80}.

\section{Background definitions and results}

\begin{Definition}
A rational matrix $G(z)\in\R(z)^{m\times m}$ is said to be {\em para-Hermitian} if $G(z) = G^*(z)$.
A para-Hermitian rational matrix $\Phi(z)\in\R(z)^{m\times m}$ is said to be a {\em spectral density} if it is positive semi-definite for all $\vartheta\in[0,2\pi)$ for which $\Phi(e^{j\vartheta})$ is defined.
A spectral density is said to be a {\em coercive} if it is positive definite in the unit circle: $\Phi(e^{j\vartheta})> 0$ for all $\vartheta\in[0,2\pi)$.
\end{Definition}

\begin{Definition}
$G(z)\in\R(z)^{m\times m}$ is said to be {\em all-pass} if 
\[
	G^*(z) G(z) = G(z) G^*(z)=I_m.
\]
Given two all-pass functions  $G_\ell(z)$ and $G_r(z)$, if $\delta_{M}(G_\ell(z))+\delta_{M}(G_r(z))=\delta_{M}(G_\ell(z)G_r(z))$, then $G_\ell(z)$ and $G_r(z)$ are said to be, respectively, {\em left all-pass divisor} and {\em right all-pass divisor}  of $G(z):=G_\ell(z)G_r(z)$.
\end{Definition}

Consider a rational spectral density $\Phi(z)\in\R(z)^{m\times m}$ of normal rank $\mathrm{rk}(\Phi)=r\leq m$. We recall that $\Phi(z)$ admits a factorization of the form \cite{BF-2016}
\[
	\Phi(z)=W(z) W^*(z),   
\]
where $W(z)\in\R^{m\times r}(z)$ is called a {\em spectral factor} of $\Phi(z)$. If $W(z)$ is such that $\delta_M(W)=\frac{1}{2}\delta_M(\Phi)$, then $W(z)$ is called a {\em minimal spectral factor} of $\Phi(z)$. We can identify four ``extremal'' minimal spectral factors of $\Phi(z)$, namely:
\begin{itemize}
  \item $W_{-}(z)$ analytic with its inverse in $\{z\in\overline{\mathbb{C}} \, :\, |z|>1\}$ ({\em minimum-phase or outer spectral factor}).

  \item $W_{+}(z)$ analytic in $\{z\in\overline{\mathbb{C}} \, :\, |z|>1\}$ with inverse analytic in $\{\,z\in\overline{\mathbb{C}} \, :\, |z|<1\,\}$.

  \item $\overline{W}_{-}(z)$ analytic in $\{z\in\overline{\mathbb{C}} \, :\, |z|<1\}$ with inverse analytic in $\{\,z\in\overline{\mathbb{C}} \, :\, |z|>1\,\}$.

  \item $\overline{W}_{+}(z)$ analytic with its inverse in  $\{z\in\overline{\mathbb{C}} \, :\, |z|<1\}$ ({\em conjugate outer spectral factor}).
\end{itemize}

These four spectral factors are connected by suitable transformations as depicted in the commutative diagram below  where an arrow indicates post-multiplication with the labelled object, e.g. $W_+(z)=W_-(z)T_1(z)$. 
\begin{center}
\begin{tikzpicture}
  \matrix (m) [matrix of math nodes,row sep=2em,column sep=2em,minimum width=2em] {
   W_{-} & & & W_{+} \\
     & & W_0 & \\
     & & & \\
     \overline{W}_{-} & & & \overline{W}_{+} \\};
  \path[-stealth]
  (m-1-1) edge[in=140,out=-20] node[fill=white] {$T_{-}$} (m-2-3)
  	  edge[in=160,out=-60] node[fill=white] {$T$} (m-4-4);
  \path[-stealth]
    (m-1-1) edge node[fill=white] {$T_1$} (m-1-4);
  \path[-stealth]
    (m-1-4) edge node[fill=white] {$T_2$} (m-4-4);   
  \path[-stealth]
    (m-2-3) edge[in=120,out=-40] node [fill=white] {$T_+$} (m-4-4);
  \path[-stealth]
  (m-1-1) edge node[fill=white]{$\overline{T}_1$} (m-4-1);
  \path[-stealth]
    (m-4-1) edge node[fill=white] {$\overline{T}_2$} (m-4-4); 
\end{tikzpicture}
\end{center}
In the next section, we will show that all minimal spectral factors are connected to $W_-(z)$ by transformations which correspond to the left all-pass divisors of the all-pass function $T(z):=W^{-L}_-(z)\overline{W}_+(z)$.
We call $T(z)$ {\em conjugate phase function} associated with the spectral density $\Phi(z)$, since it can be regarded as the conjugate version of the well-known phase function $\overline{W}^{-L}_+(z)W_{-}(z)$, which is of crucial importance in stochastic realization theory \cite{LP-15}.

Our result provides a complete parametrization of all the minimal spectral factors of a  spectral density.



\section{Main result}

\begin{Theorem}\label{thm}
  Let $\Phi(z)\in\R(z)^{m\times m}$ be a spectral density of normal rank $\rk(\Phi)=r\leq m$. Let $W_-(z)$ be the outer spectral factor of $\Phi(z)$ and $\overline{W}_+(z)$ be the conjugate outer spectral factor of $\Phi(z)$. Let $T(z):=W^{-L}_-(z)\overline{W}_+(z)$.
Let $\mathscr{W}$ be the set of minimal spectral factors of $\Phi(z)$.
Then
$$
\mathscr{W} =\{ W_-(z) T_{\ell}(z): \ T_{\ell}(z) \mbox{ is a left all-pass divisor  of } T(z)\}.
$$
\end{Theorem}

\proof

The proof is divided in three main parts:
\begin{enumerate}[I.]
\item First, we show that, without loss of generality, we can restrict the attention to spectral densities that do not have poles/zeros at infinity (set of \emph{biproper} spectral densities).
\item Then, we prove the statement of the theorem for the set of \emph{coercive} spectral densities.
\item Finally, we show how the latter result can be extended to spectral densities that have poles/zeros on the unit circle and/or are normal rank deficient.  
\end{enumerate}

{\em Part I.} As far as the first part is concerned, 
suppose that $\Phi(z)$ has a pole/zero at infinity and consider a M\"obius transformation
 $\lambda\colon \overline{\C}\to \overline{\C}$ mapping $z$ in $\lambda(z)=\frac{z-a}{1-az}$,
where $a\in\R$ is such that $|a|<1$ and $1/a$ does not coincide with a pole/zero of $\Phi(z)$.
The inverse of this map has the same structure and maps $\lambda\mapsto z(\lambda)=\frac{\lambda+a}{1+a\lambda}$. We observe that:
\begin{enumerate}
  \item $|\lambda(z)|=1$ (resp. $|\lambda(z)|>1$, $|\lambda(z)|<1$) if and only if  $|z|=1$ (resp. $|z|>1$, $|z|<1$);
  \item $[\lambda(z)]^{-1}=\lambda(z^{-1})$;
  \item $\lambda(\infty)=-1/a$ and $\lambda(z)=\infty$ if and only if $z=1/a$;
  \item If $F(z)\in\R(z)^{k\times h}$ and $G(\lambda)=F(z(\lambda))$ then $\delta_M(F(z))=\delta_M(G(\lambda))$ \cite[Chap. 4]{BGK};
  \end{enumerate}
 As a consequence
 \begin{enumerate}
 \item 
    $\Psi(\lambda):=\Phi(z(\lambda))$ is a coercive spectral density;
\item
  since $1/a$ is not a pole/zero of $\Phi(z)$ then $\Psi(\lambda)$ has no pole/zero at infinity;

  \item   the outer spectral factor $V_-(\lambda)$ (resp. conjugate outer spectral factor $\overline{V}_+(\lambda)$)
  of $\Psi(\lambda)$ is given by $V_-(\lambda)=W_-(z(\lambda))$ (resp. $\overline{V}_+(\lambda)=\overline{W}_+(z(\lambda))$);
  
 \item the conjugate phase function ${\cal T}(\lambda)$ associated with $\Psi(\lambda)$ is given by 
 ${\cal T}(\lambda)=T(z(\lambda))$ and ${\cal T}_\ell(\lambda)$ is a left all-pass divisor of ${\cal T}(\lambda)$
 if and only if ${\cal T}_\ell(\lambda)=T_{\ell}(z(\lambda))$, where $T_{\ell}(z)$ is a left all-pass divisor of $T(z)$;
 \item
 $V(\lambda)$ is a spectral factor of  $\Psi(\lambda)$ if and only if $V(\lambda)=W(z(\lambda))$ where $W(z)$ is a spectral factor of $\Phi(z)$. Moreover, in this case, $\delta_M(V(\lambda))=\delta_M(W(z))$.
 \end{enumerate}
 
Due to these facts, we can apply the argument presented in what follows to $\Psi(\lambda)$ and then transform back $\lambda(z)\mapsto z$ to recover the desired parametrization for the original spectrum $\Phi(z)$. 

{\em Part II.}  We now proceed with the second part of the proof. In this part, we will show that the theorem holds for the set of {\em coercive} spectral densities.  To this end, assume that $\Phi(z)$ is coercive. In the light of Part I, we can also assume that $\Phi(z)$ is biproper, i.e. it has no pole/zero at infinity. We first show that if $T_{\ell}(z)$ is a left all-pass divisor  of  $T(z)$ then $W(z):=W_-(z) T_{\ell}(z)$ is a minimal spectral factor of $\Phi(z)$.
To this end it is clearly sufficient to show that the McMillan degree of $W(z):=W_-(z) T_{\ell}(z)$ equals the McMillan degree of $W_-(z)$ (which, in turn, is one half of the McMillan degree of the spectral density $\Phi(z)$). To prove this fact, we start from a minimal realization of the outer spectral factor $W_-(z)$:
\beq\label{Wmeno}
W_-(z)=C(zI-A)^{-1}B+D
\eeq
and we follow five steps:
\begin{enumerate}
\item
We compute a realization of the all-pass function $T_1(z):=[W_-(z)]^{-1}W_+(z)$ in terms of the quadruple
$A,B,C,D$.
\item
We compute a realization of the all-pass function $T_2(z):=[W_+(z)]^{-1} \overline{W}_+(z)$ in terms of the quadruple
$A,B,C,D$. \item
We compute  a realization of the conjugate phase function $T(z)=T_1(z)T_2(z)$ again in terms of the quadruple
$A,B,C,D$.
\item
We use the results of \cite{Ferrante-Picci} that provide an explicit expression parametrizing the all-pass divisors of a given all-pass function; in this way, we have an expression of $T_{\ell}(z)$ in terms of the original data $A,B,C,D$ and of a free parameter.
\item We compute  the product  $W_-(z)T_{\ell}(z)$ and show that its McMillan degree equals the McMillan degree of $W_-(z)$.
\end{enumerate}

1) Let us consider a minimal realization (\ref{Wmeno}) of $W_-(z)$ and let $n$ be the McMillan degree of $W_-(z)$, i.e. the dimension of the matrix $A$. 
Let 
\beq\label{defGamma}
\Gamma:=A-BD^{-1}C
\eeq
be the zero matrix of $W_-(z)$ and recall that $\Gamma$ is non-singular and all its eigenvalues have modulus smaller than $1$. Moreover, it is worth noticing that the invertibility of $D$ in \eqref{defGamma} follows from the fact that $\Phi(z)$ is assumed to have no pole/zero at infinity.
We now show that 
\beq\label{relofT1}
T_1(z):=[W_-(z)]^{-1}W_+(z)=H_1(zI-\Gamma)^{-1}G_1+U_1
\eeq
where
\beq
\begin{array}{l}
H_1:=D^{-1}C\\
U_1:=[I+H_1X^{-1}H_1\tp]^{1/2}\\
G_1:=\Gamma X^{-1}H_1\tp U_1^{-1}
\end{array}
\eeq
and $X$ is the solution of the Stein equation
\beq\label{primalyap}
\Gamma\tp X\Gamma = X+H_1\tp H_1.
\eeq
 Before proving (\ref{relofT1}), notice that: (i)  $(A,C)$ and hence $(\Gamma, H_1)$ is observable so that $X$ is negative definite. In view of (\ref{primalyap}), this implies that $X+H_1\tp H_1$ is negative definite as well
 so that $I+H_1X^{-1}H_1\tp>0$ and hence $U_1$ and $G_1$ are well defined. (ii)
 By direct computations we get $G_1G_1\tp= \Gamma X^{-1}\Gamma\tp - X^{-1}$ so that (\ref{relofT1})
 is a minimal realization.
To prove (\ref{relofT1}), we show that: (i) the right-hand side of (\ref{relofT1}) is all-pass and
(ii) the product $W_-(z)[H_1(zI-\Gamma)^{-1}G_1+U_1]$ has a realization with the same state matrix 
$A$ of $W_-(z)$ and with zero matrix similar to $\Gamma ^{-\top}$.
As for (i), it is a matter of direct computation to show that
\beq
H_1\tp U_1=\Gamma\tp XG_1,\quad
U_1\tp U_1= I+G_1\tp XG_1.
\eeq
These conditions, together with (\ref{primalyap}), guarantee that the right-hand side of (\ref{relofT1}) is all-pass \cite[Theorem 2.1, point 3)]{Ferrante-Picci}. 
As for (ii), by taking into account that $BH_1=BD^{-1}C=A-\Gamma=(zI-\Gamma)-(zI-A)$,
we can easily see that $W_-(z)[H_1(zI-\Gamma)^{-1}G_1+U_1]=C(zI-A)^{-1}B_{+} +D_{+}$,
where 
\beq
\begin{array}{l}
B_+ :=BU_1+\Gamma X^{-1} H_1\tp U_1^{-1},\\
D_+ :=DU_1.
\end{array}
\eeq
Hence, its zero matrix is easily seen to be 
\bea
\nn
\Gamma_+ &:=& A-B_+ D_+^{-1} C\\ \nn
&=&\Gamma -\Gamma X^{-1} H_1\tp [I+H_1X^{-1}H_1\tp]^{-1} D^{-1} C\\ \nn
&=&\Gamma( X+ H_1 H_1\tp)^{-1} X=X^{-1}\Gamma ^{-\top} X.
\eea
In conclusion, $T_1(z)$ is given by the right-hand side of (\ref{relofT1}) and
\beq\label{realw+}
W_+(z)=C(zI-A)^{-1}B_{+} +D_{+}.
\eeq

2) 
We now show that 
\beq\label{realofT2}
T_2(z):=[W_+(z)]^{-1}\overline{W}_+(z)=H_2(zI-A^{-\top})^{-1}G_2+U_2
\eeq
where
$
H_2:=B_+\tp A^{-\top},\ U_2:=[I+H_2Y^{-1}H_2\tp]^{1/2},\ G_2:= A^{-\top} Y^{-1}H_2\tp U_2^{-1}
$
and $Y$ is the solution of the Stein equation
\beq\label{secondalyap}
A^{-1} YA^{-\top} = Y+H_2\tp H_2.
\eeq
Notice that from stochastic minimality of $W_+(z)$ it follows that the realization (\ref{realw+}) is minimal.
We can therefore use the same argument used in point 1) to see that: $Y$ is positive definite,
$U_2$ and $G_2$ are well defined, $U_2$ is invertible and (\ref{realofT2}) is a minimal realization.
To prove (\ref{realofT2}), we show that: (i) the right-hand side of (\ref{realofT2}) is all-pass and
(ii) the product $W_+(z)[H_2(zI-A^{-\top})^{-1}G_2+U_2]$ has a realization with  
state matrix given by 
$A^{-\top}$  and with zero matrix similar to $\Gamma ^{-\top}$.
As for (i), it is a matter of direct computation to show that
\beq
H_2\tp U_2=A^{-1} Y G_2,\quad
U_2\tp U_2= I+G_2\tp YG_2.
\eeq
These conditions, together with (\ref{secondalyap}) guarantee that the right-hand side of (\ref{realofT2}) is all-pass \cite[Theorem 2.1, point 3)]{Ferrante-Picci}.

As for (ii), by taking into account that $B_+H_2=AH_2\tp H_2=A(A^{-1} YA^{-\top} - Y) =(zI-A)Y-Y(zI-A^{-\top})$, we can easily see that
$W_+(z)[H_2(zI-A^{-\top})^{-1}G_2+U_2]=D_+ U_2 +\overline{C}_+(sI-A^{-\top})G_2 +N(z) $,
where 
\beq
\begin{array}{l}
\overline{C}_+ :=CY+D_+H_2\end{array}
\eeq
and $N(z):=C(zI-A)^{-1}(B_+ U_2-YG_2)$; it is now a matter of direct computation to show that $B_+ U_2-YG_2=0$ so that
$N(z)=0$.
The zero matrix of the product $W_+(z)[H_2(zI-A^{-\top})^{-1}G_2+U_2]$ is thus
\bea
\nn
\overline{\Gamma}_+ &:=& A^{-\top}-G_2 (D_+ U_2)^{-1} \overline{C}_+\\ \nn
&=&A^{-\top}(Y+H_2\tp H_2)^{-1} Y
\\ 
\nn
&& -  Y^{-1}(A + AH_2\tp H_2Y^{-1}) H_2\tp U_2^{-2}U_1^{-1}D^{-1}CY\\
\nn
&=&Y^{-1} A Y-  Y^{-1}AH_2\tp U_2^2 U_2^{-2}U_1^{-1}D^{-1}CY
\\ 
\nn
&=&Y^{-1} (A - AH_2\tp U_1^{-1}D^{-1}C)Y=Y^{-1} \Gamma_+ Y.
\eea
In conclusion, $T_2(z)$ is given by the right-hand side of (\ref{realofT2}).

Before proceeding to the next point, we need to establish a formula linking $X$ and $Y$.
First observe that taking the inverse of (\ref{primalyap}) and employing the Sherman-Morrison-Woodbury formula we get
\beq
X^{-1}=\Gamma X^{-1}\Gamma\tp -\Gamma X^{-1}H_1\tp (I+H_1 X^{-1}H_1\tp)^{-1}H_1  X^{-1}\Gamma\tp.
\eeq
Moreover, equation (\ref{secondalyap}) can be rewritten as 
\beq
\label{seclyap2}
Y=AYA\tp + B_+ B_+\tp.
\eeq
By direct computation, we get $B_+ B_+\tp =BB\tp +A  X^{-1}A\tp-X^{-1}$
which, plugged in (\ref{seclyap2}), gives
the  identity  
\beq
\label{terzalyap}
Z=BB\tp +A  ZA\tp,
\eeq
where
$
Z:= Y+X^{-1}.
$
Notice that by reachability of $(A,B)$, $Z$ is invertible.

3)
It is now immediate to compute the following realization of $T(z)=T_1(z)T_2(z)$:
$T(z)=\bar{\mathscr{C}}(zI-\bar{\mathscr{A}})^{-1}\bar{\mathscr{B}}+\mathscr{D}$,
where 
$\bar{\mathscr{A}}:=\bmat{cc} \Gamma &G_1H_2\\0&A^{-\top}\emat$, $\bar{\mathscr{B}}:=\bmat{cc} G_1U_2\\G_2\emat$, $\bar{\mathscr{C}}:=[H_1\mid U_1H_2]$, and $\mathscr{D}:=U_1U_2$.
By direct computation it is easy to see that 
\beq
-\Gamma X^{-1} +X^{-1} A^{-\top} + G_1H_2=0
\eeq
so that we can perform a change of basis in the state space of $T(z)$ induced by the transformation $\bar{T}= \bmat{cc} I &-X^{-1}\\0&I\emat$ in such a way that
\beq\label{minrealT}
T(z)= {\mathscr{C}}(zI- {\mathscr{A}})^{-1} {\mathscr{B}}+\mathscr{D}
\eeq
with
\beq
{\mathscr{A}}:=\bar{T}^{-1} \bar{\mathscr{A}} \bar{T}= \bmat{cc} \Gamma &0\\0&A^{-\top}\emat,
\eeq
\beq
{\mathscr{B}}:=\bar{T}^{-1} \bar{\mathscr{B}}=\bmat{c} G_1U_2+X^{-1}G_2\\G_2\emat,
\eeq
and 
\beq
{\mathscr{C}}:=\bar{\mathscr{C}} \bar{T}=[H_1\mid U_1H_2-H_1X^{-1}]=[H_1\mid B\tp A\mtp].
\eeq
Thus, is apparent that $({\mathscr{A}},{\mathscr{C}})$ is observable.
By using a dual argument it is not difficult to see that $({\mathscr{A}},{\mathscr{B}})$ is reachable so that 
(\ref{minrealT}) is a minimal realization.

Now define 
\beq
\mathscr{P}_0:=\bmat{cc} X^{-1}+X^{-1}Y^{-1}X^{-1} & X^{-1}Y^{-1} \\ Y^{-1}X^{-1}& Y^{-1} \emat
\eeq
and observe that $\mathscr{P}_0$ is invertible; in fact, 
\beq\label{eq:P0inv}
\mathscr{P}_0^{-1}=\bmat{cc} X  & -I \\ -I& Y+X^{-1} \emat=\bmat{cc} X  & -I \\ -I& Z \emat.
\eeq
By long but direct computations, we see that the following relations hold
\beq\label{lmeallpassgrande}
\left\{
\begin{array}{l}
\mathscr{A} \mathscr{P}_0 \mathscr{A}\tp -\mathscr{P}_0=\mathscr{B}\mathscr{B}\tp\\
\mathscr{A} \mathscr{P}_0 \mathscr{C}\tp= \mathscr{B}\mathscr{D}\tp\\
I+ \mathscr{C} \mathscr{P}_0 \mathscr{C}\tp= \mathscr{D}\mathscr{D}\tp
\end{array}
\right.
\eeq
Similarly, we get
\beq\label{selygrande}
\mathscr{A}\tp \mathscr{P}_0^{-1} \mathscr{A} -\mathscr{P}_0^{-1}=\mathscr{C}\tp\mathscr{C}.\\
\eeq

4) We are now in position to apply a result established in \cite[Corollary 4.1, Corollary 3.2, Remarks 3.1 and 3.2]{Ferrante-Picci}. In fact, we have a minimal realization  (\ref{minrealT}) of the all-pass function $T(z)$ and an explicit expression of the unique solution 
$\mathscr{P}_0$ of the corresponding linear matrix equation (\ref{lmeallpassgrande}).
By \cite[Corollary 4.1, Corollary 3.2, Remarks 3.1 and 3.2]{Ferrante-Picci}, 
$T_{\ell}(z)$ is a left all-pass divisor of $T(z)$ if and only if it has the form
\beq
T_{\ell}(z)= [{\mathscr{C}}(zI- {\mathscr{A}})^{-1}\mathscr{B}_{\!_\mathscr{P}}+
\mathscr{D}_{\!_\mathscr{P}} ] {\mathscr{O}}
\eeq
where
\beq\label{defbpdp}
\begin{array}{l}
\mathscr{B}_{\!_\mathscr{P}}:= \mathscr{A} \mathscr{P} \mathscr{C}\tp(I+ \mathscr{C} \mathscr{P} \mathscr{C}\tp)^{-1/2},\\
\mathscr{D}_{\!_\mathscr{P}}:=(I+ \mathscr{C} \mathscr{P} \mathscr{C}\tp)^{1/2},
\end{array}
\eeq
${\mathscr{O}}$ is an arbitrary orthogonal matrix and $\mathscr{P}$ is of the form
\beq\label{formulamitica}
\mathscr{P}=[\Pi \mathscr{P}_0^{-1} \Pi]^{+}
\eeq
with $\Pi$ (the parameter of the parametrization) varying among the orthogonal projectors on ${\mathscr{A}}$-invariant subspaces.
Notice that ${\mathscr{A}}$ is block-diagonal and its blocks $\Gamma$ and $A^{-\top}$ 
have disjoint spectra: $\sigma(\Gamma)\cap\sigma(A^{-\top})=\emptyset$.
Hence the invariant subspaces of ${\mathscr{A}}$ have the form $\im\left[\begin{array}{cc}
V_\gamma&0\\0& V_a\end{array}\right]$ where $\im(V_\gamma)$ is a $\Gamma$-invariant subspace and
$\im(V_a)$ is a $A^{-\top}$-invariant subspace.
Notice that the arbitrary orthogonal matrix ${\mathscr{O}}$ does not influence the McMillan degree of product $W_-(z)T_{\ell}(z)$ so that, without loss of generality, from now on we set  ${\mathscr{O}}=I.$

Now, let $V:=\left[\begin{array}{cc}
V_\gamma&0\\0& V_a\end{array}\right]$ be a matrix whose columns are a basis for an arbitrary 
${\mathscr{A}}$-invariant subspace.
Let $\tilde{T}:=[V\mid \tilde{W}]$, where $\tilde{W}$ is such that $\tilde{T}$ is invertible.
A change of basis on ${\mathscr{A}}$ induced by $\tilde{T}$ elicit a block-triangular structure
$\tilde{T}^{-1}{\mathscr{A}}\tilde{T}=\bmat{cc} F_{1}&\tilde{F}_{12}\\0&{F}_{2}\emat$, where, in turn, $F_{1}$ has the block-diagonal structure
\beq\label{fidiag}
F_1=\bmat{cc} \Gamma_1 &0\\0& A_1^{-\top}\emat.
\eeq
We partition now $Q:=(\tilde{T}^{-1}{\mathscr{P}_0}\tilde{T}^{-\top})^{-1}$ conformably as
$\bmat{cc} Q_{1}&Q_{12}\\Q_{12}\tp&Q_{2}\emat$. 
As shown in \cite[Lemma 4.1]{Ferrante-Picci} $Q_{1}$ is invertible so that we can set
$\bar{T}:= \bmat{cc} I&0\\-Q_{12}\tp Q_{1}^{-1}&I\emat $ and we have
$\bar{T} Q\bar{T}\tp=
\diag[ Q_{1}, Q_2-Q_{12}\tp Q_{1}^{-1}Q_{12}]$.
Therefore, by defining 
\beq\label{defofT}
T:=\tilde{T}\bar{T}\tp=[V\mid {W}]=\left[\begin{array}{ccc}
V_\gamma&0&W_1\\0& V_a&W_2\end{array}\right],
\eeq
(where $W$, partitioned conformably with $V$ in two blocks, is a new completion of $V$ to an invertible matrix)
we have
that $ {T}^{-1}{\mathscr{P}_0} {T}^{-\top}$ has the following block diagonal structure:
\beq
{T}^{-1}{\mathscr{P}_0} {T}^{-\top}=\bmat{cc} P_{1}&0\\0&P_{2}\emat.
\eeq
Moreover, the structure of ${T}^{-1}{\mathscr{A}}{T}$ is easily seen to be 
\beq\label{cbinacorsivo}
{T}^{-1}{\mathscr{A}}{T}=\bmat{cc} F_{1}&{F}_{12}\\0&{F}_{2}\emat.
\eeq
Finally by using (\ref{formulamitica}) and observing that in the new basis $\Pi$ (the orthogonal projector on $\im(V)$) is given by $\Pi=\bmat{cc} I&0\\0&0\emat$, we have
\beq
{T}^{-1}{\mathscr{P}}{T}^{-\top}=\bmat{cc} P_{1}&0\\0&0\emat.
\eeq

5) we are now ready to compute the spectral factor $W(z):=W_-(z)T_{\ell}(z)$ and show that its McMillan degree is $n$.
To this end, we first define $n_\gamma$ and $n_a$ to be the number of columns of $V_\gamma$ and $V_a$, respectively.
By direct computation, we see that
$
B\mathscr{C}=[BH_1\mid BB\tp A^{-\top}]=[BD^{-1}C \mid ZA^{-\top} -AZ]
$
or
$
B\mathscr{C}=[(zI-\Gamma)-(zI-A)  \mid  (zI-A)Z-Z(zI-A^{-\top})].
$
Moreover,
$
D\mathscr{C}=[C  \mid  DB\tp A^{-\top}].
$
It is now easy to see that
$
W(z):=W_-(z)T_{\ell}(z)=D\mathscr{D}_{\!_\mathscr{P}}+W_1(z)+W_2(z)$
where
$
W_1(z):=(CZ +DB\tp A^{-\top})(zI-A^{-\top})^{-1}[0\mid I]\mathscr{B}_{\!_\mathscr{P}}
$
and
$
W_2(z):=C(zI-A)^{-1}([I\mid -Z]\mathscr{B}_{\!_\mathscr{P}}+B\mathscr{D}_{\!_\mathscr{P}}).
$
Thus,  $\delta_M(W(z))=\delta_M(W_1(z))+\delta_M(W_2(z))$.
To compute the McMillan degree $\delta_M(W_1(z))$ of $W_1(z)$, consider the term
\begin{align*}
\nn
M&:=[0\mid I]\mathscr{B}_{\!_\mathscr{P}}
\\
&=
[0\mid I]\mathscr{A} \mathscr{P} \mathscr{C}\tp\mathscr{D}_{\!_\mathscr{P}}\inv\\
&=
[0\mid V_a\mid W_2] T\inv \mathscr{A}TT\inv \mathscr{P}T\tp T\mtp \mathscr{C}\tp\mathscr{D}_{\!_\mathscr{P}}\inv
\\
&=
[0\mid V_a\mid W_2] \bmat{cc} F_{1}&0\\0&0\emat\bmat{cc} P_{1}&0\\0&0\emat T\mtp \mathscr{C}\tp\mathscr{D}_{\!_\mathscr{P}}\inv\\
&=
[0\mid V_aA_1\mtp\mid 0] \bmat{cc} P_{1}&0\\0&0\emat T\mtp \mathscr{C}\tp\mathscr{D}_{\!_\mathscr{P}}\inv\\
&=
V_aA_1\mtp E_2,
\end{align*}
where $E_2$ is the second block rows in the partition of $\bmat{cc} P_{1}&0\\0&0\emat T\mtp \mathscr{C}\tp\mathscr{D}_{\!_\mathscr{P}}\inv$ in three block rows, consistently with the partition $[0\mid V_aA_1\mtp\mid 0]$.
From (\ref{fidiag}), (\ref{defofT}) and (\ref{cbinacorsivo}), it immediately follows that 
$A\mtp V_a =V_aA_1\mtp$
so that $(zI-A\mtp)^{-1}V_a=V_a(zI-A_1\mtp)^{-1}$.
Thus
$$
 W_1(z)  =(CZ +DB\tp A^{-\top})V_a(zI-A_1\mtp)^{-1} A_1\mtp E_2.$$
Thus, $\delta_M(W_1(z))\leq n_a=\dim(A_1)$.

To compute the McMillan degree of $W_2(z)$, we analyze the term
$N:=[I\mid -Z]\mathscr{B}_{\!_\mathscr{P}}+B\mathscr{D}_{\!_\mathscr{P}}$.
It can be rewritten as
\bea
\nn
N&=&[0\mid -I] 
\mathscr{P}_0^{-1}
\mathscr{A} \mathscr{P} \mathscr{C}\tp
\mathscr{D}_{\!_\mathscr{P}}^{-1}+B\mathscr{D}_{\!_\mathscr{P}}\\
\nn
&=&
[0\mid -I] 
(\mathscr{A}^{-\top}\mathscr{P}_0^{-1}+\mathscr{A}^{-\top}\mathscr{C}\tp\mathscr{C}) \mathscr{P} \mathscr{C}\tp
\mathscr{D}_{\!_\mathscr{P}}^{-1}+B\mathscr{D}_{\!_\mathscr{P}}
\eea
where, for the last equality we exploited (\ref{selygrande}).
By direct computation, we get $[0\mid -I] \mathscr{A}^{-\top}\mathscr{C}\tp=-B$, so that we easily obtain
$[0\mid -I] \mathscr{A}^{-\top}\mathscr{C}\tp\mathscr{C} \mathscr{P} \mathscr{C}\tp
\mathscr{D}_{\!_\mathscr{P}}^{-1}+B\mathscr{D}_{\!_\mathscr{P}}=-B(\mathscr{C} \mathscr{P} \mathscr{C}\tp
\mathscr{D}_{\!_\mathscr{P}}^{-1}-\mathscr{D}_{\!_\mathscr{P}})=
-B(\mathscr{C} \mathscr{P} \mathscr{C}\tp
-\mathscr{D}_{\!_\mathscr{P}}^2)\mathscr{D}_{\!_\mathscr{P}}^{-1}=
B\mathscr{D}_{\!_\mathscr{P}}^{-1}=[0\mid I] \mathscr{A}^{-\top}\mathscr{C}\tp\mathscr{D}_{\!_\mathscr{P}}^{-1}.$
Therefore,
$$
N=
[0\mid I] \mathscr{A}^{-\top}(I-\mathscr{P}_0^{-1}\mathscr{P}) \mathscr{C}\tp
\mathscr{D}_{\!_\mathscr{P}}^{-1}.
$$
We now use the change of basis in (\ref{defofT}) and observe that
\bea
\nn
I-\mathscr{P}_0^{-1}\mathscr{P}&=&T\mtp T\tp-T\mtp T\tp \mathscr{P}_0^{-1}T  T\inv\mathscr{P}T\mtp T\tp
\\
\nn
&=&T\mtp   \bmat{cc}0&0\\ 0 & I \emat  T\tp.
\eea
Therefore
\bea
\nn
N&=&[0\mid I]T\mtp T\tp \mathscr{A}^{-\top}T\mtp   \bmat{cc}0&0\\ 0 & I \emat  T\tp \mathscr{C}\tp
\mathscr{D}_{\!_\mathscr{P}}^{-1}
\\
\nn
&=&
[0\mid I]T\mtp \bmat{cc}F_{1}\mtp&0\\ \star & F_{2}\mtp\emat   \bmat{cc}0&0\\ 0 & I \emat  T\tp \mathscr{C}\tp
\mathscr{D}_{\!_\mathscr{P}}^{-1}\\
\nn
&=&
[0\mid I]T\mtp \bmat{cc}0 &0\\ 0 & F_{2}\mtp\emat     T\tp \mathscr{C}\tp
\mathscr{D}_{\!_\mathscr{P}}^{-1}.
\eea

Partition now $T\mtp$ conformably with $T$ as
$
T\mtp=\bmat{ccc}K_{11}\tp&K_{21}\tp&K_{31}\tp\\K_{12}\tp&K_{22}\tp&K_{32}\tp\emat
$
so that
\beq
N=
[0\mid K_{32}\tp F_{2}\mtp]     T\tp \mathscr{C}\tp
\mathscr{D}_{\!_\mathscr{P}}^{-1}
=
K_{32}\tp F_{2}\mtp[T\tp \mathscr{C}\tp
\mathscr{D}_{\!_\mathscr{P}}^{-1}]_2
\eeq
where $[T\tp \mathscr{C}\tp
\mathscr{D}_{\!_\mathscr{P}}^{-1}]_2$ denotes the second  block rows of $T\tp \mathscr{C}\tp
\mathscr{D}_{\!_\mathscr{P}}^{-1}$.
From (\ref{cbinacorsivo}) it immediately follows that $AK_{32}\tp=K_{32}\tp F_{2}\mtp$
so that $(zI-A)^{-1}K_{32}\tp=K_{32}\tp(zI-F_2\mtp)^{-1}$.
Thus
\beq\label{realw2}
 W_2(z)  =CK_{32}\tp(zI-F_2\mtp)^{-1}F_{2}\mtp[T\tp \mathscr{C}\tp
\mathscr{D}_{\!_\mathscr{P}}^{-1}]_2.
\eeq
Hence, $\delta_M(W_2(z))\leq 2n-n_a-n_\gamma=\dim(F_2)$.
To reduce this bound consider the observability matrix of realization (\ref{realw2}):
\beq
\bmat{c}
 CK_{32}\tp\\
 CK_{32}\tp F_2\mtp\\
 CK_{32}\tp(F_2\mtp)^2\\
 \vdots\emat=\bmat{c}
 C\\
 CA\\
 CA^2\\
 \vdots\emat K_{32}\tp
 \eeq
whose kernel (the un-observable subspace of the realization  (\ref{realw2})) is the kernel of 
$K_{32}\tp$.
Hence
$$
\delta_M(W_2(z))\leq 2n-n_a-n_\gamma-\dim(\ker(K_{32}\tp)).
$$
To find $\dim(\ker(K_{32}\tp))$, notice that $K_{32}\in\R^{(2n-n_a-n_\gamma)\times n}$ and from $T\inv T=I$ we immediately get $K_{32}\tp V_a=0.$ Let $\tilde{V}_a\in\R^{n\times (n-n_a)}$  be a matrix whose columns complete the columns of $V_a$ to a basis of $\R^n$ so that $[V_a\mid \tilde{V}_a]\in 
\R^{n\times n}$ is nonsingular.
We have
$
\rank(K_{32})=\rank(K_{32}[V_a\mid \tilde{V}_a])=\rank([0\mid K_{32}\tilde{V}_a])=\rank(K_{32}\tilde{V}_a)\leq n-n_a
$
because $K_{32}\tilde{V}_a \in\R^{(2n-n_a-n_\gamma)\times (n-n_a)}$.
By recalling that $\ker(K_{32}\tp)=[\im(K_{32})]^\perp$, we have
$\dim(\ker(K_{32}\tp))=2n-n_a-n_\gamma-\dim(\im(K_{32}))=2n-n_a-n_\gamma-\rank(K_{32})\geq
n-n_\gamma$.
Thus,
$
\delta_M(W_2(z))\leq 2n-n_a-n_\gamma-\dim(\ker(K_{32}\tp))\leq n-n_a.
$

In conclusion,
$
\delta_M(W(z))=\delta_M(W_1(z))+\delta_M(W_2(z))\leq n_a+ n-n_a=n, 
$ 
and hence  
$\delta_M(W(z))=n$ since $n$ is the minimal degree for a spectral factor of $\Phi(z)$.\\

We now show the opposite direction, namely that if $W_0(z)$ is a minimal spectral factor of $\Phi(z)$ then $T_-(z):=[W_-(z)]^{-1} W_0(z)$ is a left all-pass divisor of the conjugate phase function $T(z)$.
Clearly $T_-(z)$ is all-pass and by defining the all-pass function $T_+(z):=[W_0(z)]^{-1} \overline{W}_+(z)$, we have $T(z)=T_-(z) T_+(z)$
Therefore, we only need to show that $\delta_M(T_-(z)) +\delta_M(T_+(z))=\delta_M(T(z)).$ Since we have already seen that (\ref{minrealT}) is a minimal realization of $T(z)$, so that $\delta_M(T(z))=2n$,  we need to show that 
$
\delta_M(T_-(z)) +\delta_M(T_+(z))=2n.
$
But the McMillan degree of the product of two rational function is no larger than the sum of the McMillan degrees of the two factors, thus we only need to show that
\beq\label{thesivecieversa}
\delta_M(T_-(z)) +\delta_M(T_+(z))\leq 2n.
\eeq

To this aim, let us consider a minimal realization 
\beq\label{realw0}
W_0(z)=C_0(zI-A_0)^{-1} B_0+D_0
\eeq
and let $\Gamma_0:=A_0-B_0D_0^{-1}C_0$ be the corresponding zero matrix.
Notice that by the assumptions on $\Phi(z)$, $A_0,\Gamma_0$ and $D_0$ are invertible and 
none of the eigenvalues of $A_0$ and $\Gamma_0$ have unitary modulus.
We consider two different basis in the state space of $W_0(z)$: one in which
\beq\label{gamma0diag}
\Gamma_0=\bmat{cc} \Gamma_u&0\\0&\Gamma_s\emat,
\eeq
and the other
in which
$
A_0=\bmat{cc} A_u&0\\0&A_s\emat,
$
where all the eigenvalues of $\Gamma_u$ and $A_u$ have modulus larger than $1$ and
all the eigenvalues of $\Gamma_s$ and $A_s$ have modulus smaller than $1$.
Let $\gamma_1$, $\gamma_2=n-\gamma_1$, $a_1$ and $a_2=n-a_1$ be the dimensions of the matrices
$\Gamma_u,\Gamma_s,A_u$ and $A_s$, respectively.
To conclude, we show that $\delta_M(T_-(z))\leq \gamma_1+a_1$ and $\delta_M(T_+(z))\leq \gamma_2+a_2$.

Let us consider (\ref{realw0}) and the basis in which (\ref{gamma0diag}) holds.
Partition $C_0$ conformably as $C_0=[C_u\mid C_s]$.
Notice that  observability of $(A_0,C_0)$ implies observability of the pair
 $(\Gamma_0,C_0)$ and, in turn,  observability of the pair
 $(\Gamma_u,C_u)$.
Thus, equation
\beq\label{vicprimalyap}
\Gamma_u\tp X_u \Gamma_u = X_u+C_u\tp  D_0\mtp D_0^{-1}C_u
\eeq
admits a unique solution $X_u$ that is positive definite and hence invertible. 
Hence,
$
U_{-}:= [I+D_0^{-1}C_uX_{u}^{-1} C_u\tp D_0\mtp]^{1/2}$ is well defined and invertible.
Let $X_{1}:=\bmat{cc} X_u^{-1}&0\\0&0\emat$ and consider the function
\beq\label{relofT0-}
T_{1-}(z):=D_0^{-1}C_0(zI-\Gamma_0)^{-1}G_{-}+U_{-}
\eeq
where $G_{-}:=\Gamma_0 X_{1} C_0\tp D_0\mtp U_{-}^{-1}$.
Notice that $G_{-}$ can be rewritten as
$
G_{-}=\bmat{c}
\Gamma_u X_{u}^{-1} C_u\tp D_0\mtp U_{-}^{-1}\\0\emat
$
so that $T_{1-}(z)$ may also be realized as
$
T_{1-}(z):=C_{-}(zI-\Gamma_u)^{-1}B_{-}+U_{-}
$
where, $C_{-}:=D_0\inv C_u$ and $B_{-}:=\Gamma_u X_{u}^{-1} C_u\tp D_0\mtp U_{-}^{-1}$.
It is now easy to see that $T_{1-}(z)$ is all-pass. In fact, by direct computation we see that
$C_{-}\tp U_{-}=\Gamma_u\tp X_u B_{-}$ and $U_{-}\tp U_{-}=I+B_{-}\tp X_u B_{-}$ which together with (\ref{vicprimalyap}) imply that $T_{1-}(z)$ is all-pass \cite[Theorem 2.1, point 3)]{Ferrante-Picci}.
In addition, since we have derived a realization whose state matrix is $\Gamma_u$, clearly
$\delta_M(T_{1-}(z))\leq \gamma_1$.
Finally, since $U_{-}$ is invertible, $T_{1-}^{-1}(z)$ is also a proper all-pass function with McMillan degree
$\delta_M(T_{1-}^{-1}(z))=\delta_M(T_{1-}(z))\leq \gamma_1.$

We now compute $W_{0-}(z):=W_0(z)T_{1-}(z)$ which is a spectral factor of $\Phi(z)$ because 
$T_{1-}(z)$ is all-pass. By taking into account that
$B_0D_0\inv C_0=(zI-\Gamma_0)-(zI-A_0)$ a direct computation yields
\beq
W_{0-}(z)=C_0(zI-A_0)\inv B_{0-} +D_{0-},
\eeq
where $B_{0-}:=B_0U_{-}+G_{-}$ and $D_{0-}:=D_0U_{-}$.
The zero matrix $\Gamma_{-}$ of $W_{0-}(z)$ is given by
\bea
\nn
\Gamma_{-}&=&A_0 - B_{0-}D_{0-}\inv C_0=\Gamma_0-G_{-}U_{-}^{-1} D_0\inv C_0\\
\nn
&=&
\bmat{cc} \Gamma_u-\Gamma_uX_{u}^{-1} C_u\tp[D_0D_0\tp+ C_0X_{1}C_0\tp]^{-1}C_u &0\\0&\Gamma_s\emat
\eea
and, in view of (\ref{vicprimalyap}),
$$
\Gamma_{-}=\bmat{cc} X_{u}^{-1} \Gamma_u\mtp X_u &0\\0&\Gamma_s\emat.
$$
Thus all the zeros of $W_{0-}(z)$ have modulus smaller than $1$.

In conclusion, there exists a proper all-pass function $T_{1-}(z)$ with $\delta_M(T_{1-}^{-1}(z))=\delta_M(T_{1-}(z))\leq \gamma_1$, such that  $W_{0-}(z):=W_0(z)T_{1-}(z)$ is a spectral factor of $\Phi(z)$ having (i) the same state matrix of $W_0(z)$ and (ii) all its zeros inside the unit disk.

Now we consider $V_0(z):=W_{0-}^{-\ast}(z)$ which has a realization with state matrix similar to $\Gamma_-\mtp$ and zero matrix of the form 
$$
\bmat{cc} A_u\mtp &0\\0&A_s\mtp \emat,
$$
where $A_u$ has dimension $a_1$ and all its eigenvalues  have modulus larger than $1$ and
$A_s$ has dimension $n-a_1$ and all its eigenvalues  have modulus smaller than $1$.
We can apply to $V_0(z)$ the same procedure that led from $W_0(z)$ to  $W_{0-}(z)$
and we conclude that there exists a proper all-pass function $T_{2-}(z)$ with $\delta_M(T_{2-}^{-1}(z))=\delta_M(T_{2-}(z))\leq a_1,$ such that $V_0(z)T_{2-}(z)$ has (i) the same state matrix of $V_0(z)$ and (ii) all its zeros outside the unit disk.
Hence,  $W_{0-}(z)T_{2-}(z)=[V_0(z)T_{2-}(z)]^{-\ast}$ is a spectral factor of $\Phi(z)$ having
both its poles and its zeros all inside the unit disk. Hence $W_0(z)T_{1-}(z)T_{2-}(z)= W_{0-}(z)T_{2-}(z)=W_-(z)$ so that $T_-(z):=W^{-1}_-(z) W_0(z)=[T_{1-}(z)T_{2-}(z)]^{-1}$
which proves that $\delta_M(T_-(z))\leq \gamma_1+a_1$.

The same argument, this time referred to the lower blocks $A_s$ and $\Gamma_s$, now yields
$\delta_M(T_+(z))\leq \gamma_2+a_2=2n-(\gamma_1+a_1) $ and hence (\ref{thesivecieversa}).

{\em Part III.} In the last part of the proof, we will show that the result proved in Part II can be extended to the case of general spectral densities. First, by virtue of Part I, we can suppose that $\Phi(z)$ is biproper, i.e. it has no pole/zero at infinity.

In the remaining part of the proof, we will show that:
\begin{enumerate}
\item $W(z)\in\R(z)^{m\times r}$ is a minimal spectral factor of $\Phi(z)$ if and only if it can be written as $W(z)=W_{o}(z)V(z)$ where $W_{o}(z)$ is an $n\times r$ fixed factor that is biproper with zeros/poles in the unit circle and, possibly, in $z=0$, and $V(z)$ varies among the $r\times r$ biproper minimal spectral factors of a given coercive $r\times r$ spectral density $\Psi(z)$.
\item By letting $V_{-}(z)$ denote the minimum phase stable spectral factor of $\Psi(z)$, any biproper minimal spectral factor of $\Psi(z)$ can be written as $V(z)=V_{-}(z)T_{\ell}(z)$ with $T_{\ell}(z)$ being a left all pass divisor of $T(z)=W^{-L}_{-}(z)\overline{W}(z)_{+}$.
\end{enumerate} 

With reference to point 1), we first notice that the spectral density $\Phi(z)$ can be written as $\Phi(z) = F(z) D(z) C(z)$ where $F(z)\in\R[z]^{n\times r}$ and $C(z)\in\R[z]^{r\times n}$ are unimodular matrices and $D(z)\in\R(z)^{r\times r}$ is the Smith--McMillan canonical form of $\Phi(z)$ \cite{K-80}. The minimum-phase spectral factor of $\Phi(z)$ has the form (up to post-multiplication by constant orthogonal matrices) (see \cite{BF-2016}) $W_{-}(z) = F(z)\Theta(z) \Lambda_{-}(z)P_{-}(z)$, where $\Theta(z)\in\R(z)^{r\times r}$ is diagonal and has finite poles/zeros on the unit circle, $\Lambda_{-}(z)\in\R(z)^{r\times r}$ is diagonal and has all the finite (strictly) stable poles/zeros of $\Phi(z)$ in its diagonal, and $P_{-}(z)\in\R[z,z^{-1}]^{r\times r}$ is a suitable unimodular matrix. Consider the product $F(z)\Theta(z)$ which has poles/zeros in the unit circle and in $z=\infty$.  This product can be factorized as $$F(z)\Theta(z) = G_{-}(z)\Delta(z) G_{+}(z),$$ where $G_{-}(z)\in\R(z)^{n\times r}$ is biproper and of full column normal rank, $\Delta(z)\in \R(z)^{r\times r}$ is diagonal with monomials of the form $z^{\kappa_{i}}$, $\kappa_{i}\in\mathbb{Z}$, in its diagonal, and $G_{+}(z)\in\R[z,z^{-1}]^{r\times r}$ is unimodular. The previous factorization is known as a left Wiener--Hopf factorization at infinity \cite{FW-79}. Notice that $G_{-}(z)$ must have poles/zeros in the unit circle or in $z=0$ only. This follows from the fact that (i) $\Delta(z) G_{+}(z)$ can have poles/zeros in $z=0$ and in $z=\infty$ only, and (ii) the product $G_{-}(z)\Delta(z) G_{+}(z)=F(z)\Theta(z)$ has poles/zeros in the unit circle and in $z=\infty$. We define $W_{o}(z):=G_{-}(z)$, $\tilde{V}(z):=\Delta(z) G_{+}(z) \Lambda_{-}(z) P_{-}(z)$, and $\Psi(z) := \tilde{V}(z)\tilde{V}^{*}(z)$.
Notice that $\Psi(z)$ has no zeros/poles in the unit circle and it has full normal rank, that is, it is coercive. Let  $\{\beta_{i}\}_{i=1}^{q}$ denote the poles in the unit circle of $\Phi(z)$ and $\{\alpha_{i}\}_{i=1}^{p}$ denote the remaining poles of $\Phi(z)$. Note that, by construction, it holds 
\begin{align}\label{eq:themdeg}
&\delta(\Phi;\alpha_{i})=\delta(\Psi;\alpha_{i}),\ i=1,\dots,p, \\
&\delta(W_{o};\beta_{i})=\frac{1}{2}\delta(\Phi;\beta_{i}),\ i=1,\dots,q.
\end{align}
Observe also that, besides the poles $\{\alpha_{i}\}_{i=1}^{p}$, $\Psi(z)$ can possess additional poles only in $z=0$ and $z=\infty$.

Now, let $V(z)$ be any biproper minimal spectral factor of $\Psi(z)$, and consider $$W(z) := W_{o}(z) V(z).$$
Since $V(z)$ is taken to be minimal and biproper then, all its non-zero poles are in $\{\alpha_{i}\}_{i=1}^{p}$ and they satisfy $\sum_{i=1}^{p} \delta(V;\alpha_{i})=\frac{1}{2}\sum_{i=1}^{p} \delta(\Psi;\alpha_{i})$ \cite[Sec.~9.1]{BGKR}. Now, notice that $W(z)$ is a spectral factor of $\Phi(z)$ that is again biproper, since $W_{o}(z)$ and $V(z)$ are so. Moreover $W(z)$ has no pole/zero in $z=0$. This follows from the fact that (i) $\Phi(z)$ is biproper by assumption, and (ii) $W^{*}(z)$ has no pole/zero in $z=0$.
In view of the previous observations and of \eqref{eq:themdeg},
\begin{align*}
	\delta_{M}(W) 
	&= \sum_{i=1}^{p} \delta(V;\alpha_{i}) + \sum_{i=1}^{q} \delta(W_{o};\beta_{i}) \\
	& = \frac{1}{2}\sum_{i=1}^{p} \delta(\Psi;\alpha_{i}) + \frac{1}{2}\sum_{i=1}^{q} \delta(\Phi;\beta_{i}) \\
	& = \frac{1}{2}\sum_{i=1}^{p} \delta(\Phi;\alpha_{i}) + \frac{1}{2}\sum_{i=1}^{q} \delta(\Phi;\beta_{i}) = \frac{\delta_{M}(\Phi)}{2},
\end{align*}
i.e.,  $W(z)$ is a minimal spectral factor of $\Phi(z)$.

Conversely, let $W(z)$ be any minimal spectral factor of $\Phi(z)$ and let $V_{-}(z)$ be the stable minimum-phase spectral factor of $\Psi(z)$. As shown before, $W_{o}(z)V_{-}(z)$ is a minimal spectral factor of $\Phi(z)$, which in this case coincides with the stable minimum-phase one, that is $W_{-}(z)=W_{o}(z)V_{-}(z)$. It holds $$W(z) = W_{-}(z)U(z) = W_{o}(z)V_{-}(z)U(z),$$
for a suitable all pass matrix $U(z)\in\R(z)^{r\times r}$. Now observe that:
\begin{enumerate}[(i)]
\item $W$ has poles only in $\{\alpha_{i}\}_{i=1}^{p}\cup \{\beta_{i}\}_{i=1}^{q}$ since it is minimal,
\item $U$ cannot have poles/zeros in the unit circle since $U(e^{j\vartheta})U^{*}(e^{j\theta})=I_{r}$ for every $\theta\in[0,2\pi)$,
\item $W_{o}$ has poles/zeros in the unit circle and in $z=0$ only.
\end{enumerate} 
These three facts together imply that the non-zero poles of $V_{-}(z)U(z)$ belong to $\{\alpha_{i}\}_{i=1}^{p}$ and they satisfy $\delta(W;\alpha_{i}) = \delta(V_{-}U;\alpha_{i})$, $i=1,\dots,p$.
Moreover, we have that
\begin{align}\label{eq:theconv}
	&\sum_{i=1}^{p} \delta(V_{-}U;\alpha_{i})  = \sum_{i=1}^{p}\delta(W;\alpha_{i}) = \frac{1}{2}\sum_{i=1}^{p}\delta(\Phi;\alpha_{i}) \overset{\eqref{eq:themdeg}}{=} \frac{1}{2}\sum_{i=1}^{p}\delta(\Psi;\alpha_{i}).
\end{align}
Finally, we notice that $U(z)$ cannot have poles/zeros in $z=0$ and $z=\infty$ (i.e. $U(z)$ must be biproper), otherwise $W(z)=W_{-}(z)U(z)$ would have a pole in $z=0$ or $z=\infty$ and consequently, in view of the biproperness of $\Phi(z)$, it would not be minimal. This implies that $\delta(V_{-}U;0)=\delta(\Psi;0)$.
The latter observation together with \eqref{eq:theconv} yields 
\begin{align*}
	\delta_{M}(V_{-}U)&=\sum_{i=1}^{p} \delta(V_{-}U;\alpha_{i})  +  \delta(V_{-}U;0) \\
	&= \frac{1}{2}\sum_{i=1}^{p}\delta(\Psi;\alpha_{i}) + \delta(\Psi;0) =\frac{1}{2}\delta_{M}(\Psi),
\end{align*}
i.e. the product $V_{-}(z)U(z)$ is a biproper minimal spectral factor of $\Psi(z)$. \\

We now address point 2) We first notice that 
\begin{align}\label{eq:T}
	T(z)&=W_{-}^{-L}(z)\overline{W}_{+}(z) = V_{-}^{-1}(z)W_{o}^{-L}(z)W_{o}(z) \overline{V}_{+,0}(z) \nonumber \\
	&= V_{-}^{-1}(z) \overline{V}_{+,0}(z),
\end{align}
where $\overline{V}_{+,0}(z)$ denotes the minimal biproper spectral factor of $\Psi(z)$ having unstable poles/zeros with the only exception for those in $z=0$. Since $\Psi(z)$ is coercive, we can apply Theorem \ref{thm} and conclude that any minimal spectral factor  $V_{\ell}(z)$ of $\Psi(z)$ (and in particular the biproper ones) can be written in the form $$V_{\ell}(z) = V_{-}(z) T_{\ell}(z),$$ 
with $T_{\ell}(z)$ being a left all pass divisor of $\overline{T}(z):=V_{-}^{-1}(z) \overline{V}_{+}(z)$, i.e. $\overline{T}(z) = T_{\ell}(z) T_{\ell,r}(z)$ with $\delta_{M}(\overline{T})=\delta_{M}(T_{\ell}) + \delta_{M}(T_{\ell,r})$, where in this case $\overline{V}_{+}(z)$ denotes the unstable maximum-phase spectral factor of $\Psi(z)$. In particular, when applied to $\overline{V}_{+,0}(z)$ the latter result reads as $$\overline{V}_{+,0}(z) = V_{-}(z) T(z),$$
where $T(z)$, as defined in \eqref{eq:T}, must be a left all pass divisor of $\overline{T}(z)$, i.e. $\overline{T}(z) = T(z) T_{r}(z)$ with $\delta_{M}(\overline{T})=\delta_{M}(T) + \delta_{M}(T_{r})$. We claim that $T_{r}(z)$ can have poles/zeros in $z=0$ and $z=\infty$ only. To prove this fact, we decompose $\Psi(z)$ in its Smith--McMillan form $\Psi(z)=\tilde{F}(z)\tilde{D}(z)\tilde{C}(z)$, where $\tilde{F}(z)$ and $\tilde{C}(z)$ are $r\times r$ unimodular matrices and $\tilde{D}(z)$ is the Smith--McMillan canonical form of $\Psi(z)$. In \cite{BF-2016} it is shown that (up to post-multiplication by orthogonal matrices) $\overline{V}_{+,0}(z)$ and $\overline{V}_{+}(z)$ have the form   $\overline{V}_{+,0}(z)=\tilde{F}(z)\bar{\Lambda}_{+,0}(z)P_{+,0}(z), \overline{V}_{+}(z)= \tilde{F}(z)\bar{\Lambda}_{+}(z)P_{+}(z)$, where $\bar{\Lambda}_{+,0}(z)$ and $\bar{\Lambda}_{+}(z)$ are the diagonal parts of $\tilde{D}(z)$ containing the zeros/poles in $\{z\in\C: |z|>1 \cup z=0\}$ and in $\{z\in\C: |z|>1\}$,  respectively, and $P_{+}(z)$ and $P_{+,0}(z)$ are unimodular matrices. In the light of this fact, we have that 
\[
	T_{r}(z) = \overline{V}_{+}^{-1}(z)\overline{V}_{+,0}(z) = P_{+}(z)^{-1}\bar{\Lambda}_{+}^{-1}(z) \bar{\Lambda}_{+,0}(z)P_{+,0}(z),
\]
is a matrix whose poles are in $z=0$ or $z=\infty$ only.

Now, for every biproper minimal spectral factor $V_{\ell}(z)$ of $\Psi(z)$, we have 
$$
	V_{\ell}(z) = V_{-}(z) T_{\ell}(z), \ T(z)  = T_{\ell}(z) T_{r,0,\ell}(z),
$$
where  $T_{r,0,\ell}(z) := T_{\ell,r}(z)T_{r}^{-1}(z)$. Since $T_{\ell}(z)$ and $T(z)$ have no pole/zero in $z=0$ and $z=\infty$, it follows that $T_{r,0,\ell}(z)$ must be biproper, so that all the zeros/poles in $z=0$ and $z=\infty$ must cancel out in the product $T_{\ell,r}(z)T_{r}^{-1}(z)$. 
Let $\{\gamma\}_{i=1}^{t}$ be the poles of $\overline{T}(z)$ different from $0$ and $\infty$. In view of the minimality of the factorization $\overline{T}(z)=T_{\ell}(z) T_{\ell,r}(z)$, we have $\delta(\overline{T};\gamma_{i}) = \delta(T_{\ell};\gamma_{i}) + \delta(T_{\ell,r};\gamma_{i})$, $i=1,\dots,t,$ \cite[Sec.~9.1]{BGKR}.
Since $T_{r}^{-1}(z)$ has poles only in $z=0$ and in $z=\infty$ and $T(z)$, $T_{r,0,\ell}(z)$ are biproper, it follows that (i) $T(z)=\overline{T}(z)T_{r}^{-1}(z)$ has the same poles (and polar degrees) of $\overline{T}(z)$ except for those  in $0$ and $\infty$, and (ii) $T_{r,0,\ell}(z) := T_{\ell,r}(z)T_{r}^{-1}(z)$ has the same poles (and polar degrees) of $T_{\ell,r}(z)$ except for those in $0$ and $\infty$. This  implies
\begin{align*}
	\delta_{M}(T) &=
	 \sum_{i=1}^{t} \delta(\overline{T};\gamma_{i}) = \sum_{i=1}^{t} \delta(T_{\ell};\gamma_{i}) + \delta(T_{\ell,r};\gamma_{i}) \\
	&= \sum_{i=1}^{t} \delta(T_{\ell};\gamma_{i}) + \delta(T_{\ell,r,0};\gamma_{i}) =\delta_{M}(T_{\ell}) + \delta_{M}(T_{\ell,r,0}).
\end{align*}

Therefore any minimal biproper spectral factor of $\Psi(z)$ can be written as $V_{\ell}(z) = V_{-}(z) T_{\ell}(z)$ with $T_{\ell}(z)$ being a left all pass divisor of $T(z)$. Eventually, by virtue of the one-to-one relation between the biproper minimal spectral factors of $\Psi(z)$ and the minimal spectral factors of $\Phi(z)$ the latter result applies to the minimal spectral factors of $\Phi(z)$ as well.
\qed

 \section{A numerical example}
 
In this section, we apply our main result to a concrete example arising from stochastic realization theory. To this end, consider a zero-mean purely nondeterministic second-order stationary process $\{y(t)\}$ 
whose spectral density is
\[
	\Phi(z)=\frac{1}{z^2-\frac{5}{2}z+1}\begin{bmatrix}3z^{2}+\frac{17}{8}z+3 & 0 \\ 0 & \frac{2}{3}z^{2}-\frac{20}{9}z+\frac{2}{3}\end{bmatrix}.
\]
We want to compute {\em all} the minimum ``complexity'' (i.e., with minimal McMillan degree) dynamical models for the process $y$.
One possible model is the following {\em minimum phase} model that can be computed by standard procedures:
\begin{align*}
\begin{cases}
x(t+1) = Ax(t)+Bu(t),\\
y(t) = Cx(t) + Du(t),
\end{cases}
\end{align*}
where
\[
	A=\frac{1}{2}I_{2}, \ B=I_{2}, \ C=\begin{bmatrix}1/4& 0 \\ 0 & 1/6\end{bmatrix}, \ D=I_{2},
\]
 and with $\{u(t)\}$ being a white noise process. 
This is, however, just one possible choice: to obtain all the models of minimum complexity we can apply our result as follows.

Building on the second part of the proof of Theorem~\ref{thm}, we compute the conjugate phase function $T(z)$ that admits the following minimal state space realization
\[
	T(z)=\mathscr{C}(zI_{4}-\mathscr{A})^{-1}\mathscr{B}+\mathscr{D}.
\]
where 
\begin{align*}
	 \mathscr{A}&=\begin{bmatrix}1/4 & 0 & 0 & 0 \\ 0 & 1/3 & 0 & 0 \\ 0 & 0 & 2 & 0 \\ 0 & 0 & 0 & 2\end{bmatrix}, \ \mathscr{B}=\begin{bmatrix}-15/14 & 0 \\ 0 & -16/15 \\ -3/7 & 0 \\ 0 & -3/10\end{bmatrix},\\
	\mathscr{C}&=\begin{bmatrix}1/4 & 0 & 2 & 0\\ 0 & 1/6 & 0 & 2\end{bmatrix}, \ \mathscr{D}=\begin{bmatrix}1/2 & 0\\ 0 & 2/3\end{bmatrix}.
\end{align*}

From \cite[Corollary 4.1]{Ferrante-Picci}, we know that there is a one-to-one correspondence between left all-pass divisors of $T(z)$ and invariant subspaces of $\mathscr{A}$. The invariant subspaces of $\mathscr{A}$ can be classified in four ``classes'', namely
\[
	\{0\}\oplus \mathscr{V}, \ \mathrm{span}\begin{bmatrix}1 \\ 0\end{bmatrix}\oplus \mathscr{V}, \ \mathrm{span}\begin{bmatrix}0\\ 1\end{bmatrix}\oplus \mathscr{V}, \ \R^{2}\oplus \mathscr{V},
\]
where $\mathscr{V}$ is any subspace of $\R^{2}$ and the symbol $\oplus$ denotes direct sum of subspaces. In view of Theorem~\ref{thm}, to each element of these classes there corresponds 
an (essentially unique) minimal spectral factor of $\Phi(z)$. Notice that, in this case, each class generates an \emph{infinite} number of (essentially unique) minimal spectral factors. We now explicitly compute the minimal spectral factors belonging to each class. 

Here, we focus only on the first class, the others being similar. Thus, we consider the class of $\mathscr{A}$-invariant subspaces
\[
	\{0\}\oplus \mathscr{V}.
\]
We can further divide this class into three sub-classes:
\begin{enumerate}
\item $\{0\}\oplus \{0\}$.
\item $\{0\}\oplus \R^{2}$.
\item $\mathscr{W}_{\theta}:=\{0\}\oplus \mathrm{span}\begin{bmatrix}\cos(\theta) \\ \sin(\theta)\end{bmatrix}$, $\theta\in[0,\pi)$, where we let $\cos(\theta):=\mathrm{cos}(\theta)$ and $\sin(\theta):=\mathrm{sin}(\theta)$ to simplify the notation.
\end{enumerate}
In the first case, $T_{\ell}=I_{2}$ and the corresponding spectral factor is the minimum-phase one, namely $W_{-}(z)$. In the second and third case, we first compute  $\mathscr{P}_{0}^{-1}$ according to the expression \eqref{eq:P0inv} derived in the proof
\[
	\mathscr{P}_{0}^{-1}=\begin{bmatrix}1/15 & 0 & -1 & 0\\ 0 & 1/32 & 0 & -1\\ -1 & 0 & 4/3 & 0 \\ 0 & -1 & 0 & 4/3\end{bmatrix}.
\]
Then, for the second case, we have 
\[
	\mathscr{P} =\left[\Pi_{2}\mathscr{P}_{0}^{-1}\Pi_{2}\right]^{+}=\begin{bmatrix}0 & 0 & 0 & 0\\ 0 & 0 & 0 & 0\\ 0 & 0 & 3/4 & 0 \\ 0 & 0 & 0 & 3/4\end{bmatrix},
\]
where $\Pi_{2}$ denotes the orthogonal projection onto $\{0\}\oplus \R^{2}$.
This yields the (essentially unique) all-pass divisor
\[
	\overline{T}_{1}(z) = \mathscr{C}_{1}(zI_{2}-\mathscr{A}_{1})^{-1}\mathscr{B}_{1} +\mathscr{D}_{1},
\]
with 
\begin{align*}
	\mathscr{A}_{1}=2I_{2},\ \mathscr{B}_{1}=\frac{3}{2}I_{2}, \ \mathscr{C}_{1}=2I_{2}, \  \mathscr{D}_{1}=2I_{2}.
\end{align*}
To this all-pass divisor there corresponds the (essentially unique) unstable minimum-phase spectral factor 
\[
	\overline{W}_{-}(z)=W_{-}(z)\overline{T}_{1}(z) = \bar{C}(zI_{2}-\bar{A})^{-1}\bar{B} +\bar{D},
\]
where
\begin{align*}
	\bar{A}=2I_{2}, \ \bar{B}=\begin{bmatrix}-{4}/{5} & {8}/{5} \\ -{8}/{5} & -{4}/{5} \end{bmatrix}, \ \bar{C}=\begin{bmatrix}-{7}/{8} & -{7}/{4} \\ {5}/{3} & -{5}/{6}\end{bmatrix}, \ \bar{D}=2I_{2}. 
\end{align*}
In the third case, we have 
\[
	\mathscr{P}_{\theta} =\left[\Pi_{\mathscr{W}_{\theta}}\mathscr{P}_{0}^{-1}\Pi_{\mathscr{W}_{\theta}}\right]^{+}=\begin{bmatrix}0 & 0 & 0 & 0\\ 0 & 0 & 0 & 0\\ 0 & 0 & \frac{3}{4}c^{2}(\theta) & \frac{3}{4}\cos(\theta)\sin(\theta) \\ 0 & 0 & \frac{3}{4}\cos(\theta)\sin(\theta) & \frac{3}{4}\sin^{2}(\theta)\end{bmatrix},
\]
where $\Pi_{\mathscr{W}_{\theta}}$ denotes the orthogonal projection onto $\mathscr{W}_{\theta}$. This yields an entire family of (essentially unique) all-pass divisors parametrized by $\theta$
\[
	T_{\theta}(z) = \mathscr{C}_{1,\theta}(z-2)^{-1}\mathscr{B}_{1,\theta} +\mathscr{D}_{1,\theta},
\]
with 
\begin{align*}
	\mathscr{B}_{1,\theta}&=\begin{bmatrix}\cos(\theta)\\ \sin(\theta) \end{bmatrix}^{\top}, \ \mathscr{C}_{1,\theta}=3\begin{bmatrix}\cos(\theta) \\ \sin(\theta) \end{bmatrix},  \\ \mathscr{D}_{1,\theta}&=\begin{bmatrix}1+\cos^{2}(\theta) & \cos(\theta)\sin(\theta)\\ \cos(\theta)\sin(\theta) & 1+\sin^{2}(\theta)\end{bmatrix}.
\end{align*}
As before, to this family of all-pass divisors there corresponds the family of (essentially unique) minimal spectral factors 
\[
	W_{\theta}(z)=W_{-}(z)T_{\theta}(z)=C_{\theta}(zI_{2}-A_{\theta})^{-1}B_{\theta}+D_{\theta}, \ \theta\in[0,\pi),
\]
where the matrices $A_{\theta}$, $B_{\theta}$, $C_{\theta}$, and $D_{\theta}$ are defined as
\begin{align*}
	A_{\theta}&=\begin{bmatrix} \frac{8 \sin^2(\theta)+2}{\sin^2(\theta)+4} &-\frac{3 \sqrt{5} \sin(2 \theta)}{\cos(2 \theta)-9} \\  -\frac{3 \sqrt{5} \sin(2 \theta)}{\cos(2 \theta)-9} & -\frac{11 \cos(2 \theta)+21}{2 (\cos(2 \theta)-9)} \end{bmatrix}, \\ B_{\theta}&=\begin{bmatrix} \frac{3 \cos(\theta) \sin(\theta)}{\sqrt{\sin^2(\theta)+4}} &  \frac{3 \sin^2(\theta)+2}{\sqrt{\sin^2(\theta)+4}} \\ \frac{\sqrt{\frac{5}{2}} (\cos(2 \theta)+3)}{\sqrt{9-\cos(2 \theta)}} & \frac{\sqrt{\frac{5}{2}} \cos(\theta)\sin(\theta) \sqrt{9-\cos(2 \theta)} }{\sin^2(\theta)+4}\end{bmatrix}, \\ 
	 C_{\theta}&=\begin{bmatrix} \frac{3 \cos(\theta) \sin(\theta)}{\sqrt{\sin^2(\theta)+4}} & \frac{\sqrt{9-\cos(2 \theta)} (23 \cos(2 \theta)+33)}{8 \sqrt{10} \left(\sin^2(\theta)+4\right)} \\  \frac{9 \sin^2(\theta)+1}{3 \sqrt{\sin^2(\theta)+4}} &\frac{7 \sqrt{\frac{5}{2}} \cos(\theta)\sin(\theta) \sqrt{9-\cos(2 \theta)} }{6 \left(\sin^2(\theta)+4\right)}\end{bmatrix}, \\ D_{\theta}&=\begin{bmatrix}1+\cos^{2}(\theta) & \cos(\theta)\sin(\theta)\\ \cos(\theta)\sin(\theta) & 1+\sin^{2}(\theta)\end{bmatrix}. 
\end{align*}
Similarly, we can obtain a parametrization of each one of the other three classes of spectral factors and in such a way we get four classes of spectral factors accounting for {\em all} the models of minimal complexity for $y$.
Notice, that each class contains infinitely many (essentially different) spectral factors so that we have parametrized all the {\em infinitely many} models of minimal complexity for $y$.

\section{Conclusions}

In this paper we provide a parametrization of the set of minimal spectral factors of a discrete-time spectral density in terms of the all-pass divisors of an all-pass function (the conjugate phase function). Remarkably, our main theorem applies to general spectral densities and gives an answer to a conjecture of \cite{BF-2016}.
Moreover, this result is particularly interesting in the light of the recent work \cite{Ferrante-Picci}. In fact, in \cite{Ferrante-Picci} the set of all-pass divisors of a given all-pass function is parametrized both algebraically (in terms of solutions of a certain Algebraic Riccati Equation) and geometrically (in terms of invariant subspaces of a certain matrix). These parametrizations are therefore inherited by the set of spectral factors of a coercive spectral density. 
A final comment on the fact that our theory applies to classical spectral factorization which, beside stochastic  realization, is relevant for ${\mathcal H}_2$ problems; however, as shown e.g. in \cite{Colaneri-F-siam}, spectral factorization can also be regarded as an intermediate step to compute $J$-spectral factors associated to ${\mathcal H}_\infty$ problems.


\begin{thebibliography}{100}

\bibitem{AM-79}
B.~D.~O.~Anderson and J.~B.~Moore.
\newblock {\em Optimal Filtering.}
\newblock Prentice-Hall, Englewood Cliffs, NJ, 1979.

\bibitem{BGK}
H.~Bart, I.~Gohberg, and M.~Kaashoek. 
\newblock {\em Minimal factorization of matrix and operator functions.}
\newblock Birkh\"auser Verlag, Basel, 1979.

\bibitem{BGKR} H.~Bart, I.~Gohberg, M.~Kaashoek, and A.~C.~Ran. 
\newblock {\em Factorization of matrix and operator functions: the state space method}
\newblock  Springer, 2007.

\bibitem{BF-2016}
G.~Baggio and A.~Ferrante. 
\newblock On the factorization of rational discrete-time spectral densities.
\newblock {\em IEEE Trans. Automat. Contr.}, 
\newblock 61(4):969--981, 2016.

\bibitem{BF-min-2016}
G.~Baggio and A.~Ferrante. 
\newblock On minimal spectral factors with zeroes and poles lying on prescribed regions.
\newblock {\em IEEE Trans. Automat. Contr.}, 
\newblock 61(8):2251--2255, 2016.


\bibitem{C-93}
D.~J.~Clements.
\newblock Rational spectral factorization using state-space methods,
\newblock {\em Systems \& control letters},
\newblock 20(5):335--343, 1993.

\bibitem{Colaneri-F-siam}
P.~Colaneri and A.~Ferrante.
\newblock {Algebraic Riccati Equation and $J$-Spectral Factorization for $\mathcal{H}_{\infty}$ Filtering and Deconvolution}.
\newblock {\em SIAM J. Contr. and Opt.},
\newblock Vol. 45(1):123--145, 2006.


\bibitem{F-94}
A.~Ferrante. 
\newblock A parametrization of minimal stochastic realizations.
\newblock {\em IEEE Trans. Automat. Contr.}, 
\newblock 39(10):2122--2126, 1994.

\bibitem{F-97}
A.~Ferrante. 
\newblock A Parametrization of the Minimal Square Spectral Factors of a Nonrational Spectral Density. 
\newblock {\em J. Math. Systems, Estimation, and Control}, 
\newblock 7(2):197--226, 1997.

\bibitem{F-05}
A.~Ferrante.
\newblock Minimal representations of continuous-time processes having spectral density with zeros in the extended imaginary axis. 
\newblock {\em Systems \& control Letters}, 
\newblock 54(5):511--520, 2005.

\bibitem{FMP-93}
A.~Ferrante, G.~Michaletzky, and M.~Pavon. 
\newblock Parametrization of all minimal square spectral factors. 
\newblock {\em System \& Control Letters}, 
\newblock21:249--254, 1993.

\bibitem{Ferrante-Ntog-Automatica-13}
A.~Ferrante and L.~Ntogramatzidis.
\newblock {The Generalised Discrete Algebraic Riccati Equation in Linear-Quadratic Optimal Control}.
\newblock {\em Automatica},
\newblock 49:471--478, 2013.

\bibitem{Ferrante-Picci}
A.~Ferrante and G.~Picci.
\newblock {Representation and Factorization of Discrete-Time Rational All-Pass Functions}.
\newblock {\em IEEE Trans. Automat. Contr.},
\newblock 62(7):3262--3276, 2017.

\bibitem{FP-82}
L.~Finesso and G.~Picci. 
\newblock A characterization of minimal square spectral factors. 
\newblock {\em IEEE Trans. Automat. Contr.}, 
\newblock 27(1):122--127, 1982.

\bibitem{Furhmann-95}
P.~A.~Fuhrmann. 
\newblock On the characterization and parametrization of minimal spectral factors.
\newblock {\em J. Math. Systems, Estimation, and Control,}
\newblock 5:383--444, 1995.

\bibitem{FG-98}
P.~A.~Fuhrmann and A.~Gombani. 
\newblock On a Hardy space approach to the analysis of spectral factors.
\newblock {\em Int. J. Control,}
\newblock 71(2):277--357, 1998.

\bibitem{FW-79}
P.~A.~Fuhrmann and J.~C.~Willems.
\newblock Factorization indices at infinity for rational matrix functions.
\newblock {\em Integral equations and operator theory}
\newblock 2.3 (1979): 287-301.


\bibitem{K-80}
T.~Kailath.
\newblock {\em Linear systems}.
\newblock Vol. 156. Englewood Cliffs, NJ: Prentice-Hall, 1980.

\bibitem{Lindquist-P-79-siam}
A.~Lindquist and G.~Picci.
\newblock On the stochastic realization problem.
\newblock {\em SIAM J. Contr. and Opt.}, 17(3):365--389, 1979.

\bibitem{Lindquist-P-85-siam}
A.~Lindquist and G.~Picci.
\newblock Realization theory for multivariate stationary gaussian processes.
\newblock {\em SIAM J. Contr. and Opt.}, 23(6):809--857, 1985.

\bibitem{Lindquist-P-91-jmsec}
A.~Lindquist and G.~Picci.
\newblock A geometric approach to modeling and estimation of linear stochastic systems.
\newblock {\em J. Math. Systems, Estimation, and Control}, 1:241--333, 1991.

\bibitem{LP-15}
A.~Lindquist and G.~Picci. 
\newblock {\em Linear Stochastic Systems: A Geometric Approach to Modeling, Estimation and Identification.} 
\newblock Series in Contemporary Mathematics, Springer, 2015.


\bibitem{O-05}
C.~Oar\u{a}.
\newblock Constructive solutions to spectral and inner-outer factorizations with respect to the disk,
\newblock {\em Automatica},
\newblock 41(11):1855--1866, 2005.

\bibitem{OA-11}
C.~Oar\u{a} and R.~Andrei.
\newblock Computation of the general $(J, J^{\prime})$-lossless factorization,
\newblock {\em IEEE Trans. Automat. Contr.},
\newblock 49(7):710-717, 2013.

\bibitem{OM-13}
C.~Oar\u{a} and R.~Marinic\u{a}.
\newblock $J$ factorizations of a general discrete-time system,
\newblock {\em Automatica},
\newblock 49(7):2221-2228, 2013.

\bibitem{P-93}
M.~Pavon. 
\newblock On the parametrization of non-square spectral factors, 
\newblock in: U. Helmke, R. Mennicken, J. Sauer (Eds.), {\em Systems and Networks: Mathematical Theory and Applications}, vol. II, Mathematical Research, vol. 79, Akademie Verlag, Berlin, 1994, pp. 413--416.

\bibitem{P-Ran-01}
M.~A.~Petersen and A.~C.~M.~Ran. 
\newblock Minimal square spectral factors via triples. 
\newblock {\em SIAM J. Matrix Anal. Appl.},
\newblock 22(4):1222--1244, 2001.

\bibitem{P-Ran-02}
M.~A.~Petersen and A.~C.~M.~Ran. 
\newblock Minimal nonsquare spectral factors.
\newblock {\em Linear Algebra Appl.}, 
\newblock 351 (2002): 553--565.

\bibitem{P-Ran-02b}
M.~A.~Petersen and A.~C.~M.~Ran. 
\newblock Nonsquare spectral factors via factorizations of a unitary function.
\newblock {\em Linear Algebra Appl.}, 
\newblock 351 (2002): 567--583.

\bibitem{Picci-P-94}
G.~Picci and S.~Pinzoni.
\newblock Acausal models and balanced realizations of stationary processes.
\newblock {\em Linear Algebra Appl.}, 205-206:997--1043, 1994.

\bibitem{Ran-95}
A.~C.~M.~Ran. 
\newblock Minimal square spectral factors. 
\newblock {\em Systems \& control letters},
\newblock 24(5):307--316, 1995.

\bibitem{Willems-1971}
J.~C.~Willems.
\newblock Least squares stationary optimal control and the algebraic Riccati equation.
\newblock {\em IEEE Trans. Automat. Contr.}, 16(6): 621--634, 1971.

\bibitem{Youla-1961}
D.~C.~Youla. 
\newblock On the factorization of rational matrices. 
\newblock {\em  IRE Trans. Information Theory}, 7(3):172--189, 1961.

\bibitem{ZDG-96}
K.~Zhou, J.~C.~Doyle, and K.~Glover. 
\newblock {\em Robust and optimal control.} 
\newblock Vol. 40. New Jersey: Prentice hall, 1996.

\end{thebibliography}
\end{document}